\documentclass{amsart}
\usepackage{amssymb, amsmath, amsfonts, amscd}

\input xypic

\theoremstyle{plain}
\newtheorem{theorem}{Theorem}[section]
\newtheorem{corollary}[theorem]{Corollary}
\newtheorem{lemma}[theorem]{Lemma}
\newtheorem{proposition}[theorem]{Proposition}

\newtheorem{example}[theorem]{Example}

\theoremstyle{definition}
\newtheorem{definition}[theorem]{Definition}

\theoremstyle{remark}

\numberwithin{equation}{theorem}

\newcommand{\I}{\mathcal{I}}
\newcommand{\J}{\mathcal{J}}

\renewcommand{\L}{\mathcal{L}}
\newcommand{\E}{\mathcal{E}}

\renewcommand{\O}{\mathcal{O} }

\renewcommand{\P}{\mathbf{P} }

\renewcommand{\H}{\operatorname{H} }

\newcommand{\Pic}{\operatorname{Pic} }

\newcommand{\Sym}{\operatorname{Sym}}

\newcommand{\Ext}{\operatorname{Ext} }

\newcommand{\Z}{\mathbf{Z} }

\newcommand{\K}{\operatorname{K}}
\newcommand{\Spec}{\operatorname{Spec}}

\begin{document}

\title{Jet bundles on projective space: New examples}

\author{Helge Maakestad}
\address{Hoegskolen i Bergen}

\email{\text{h\_maakestad@hotmail.com} }


\date{Spring 2012}

\begin{abstract} In a previous paper on this subject the structure of the sheaf of principal parts as left and right module over the structure sheaf 
on the projective line was studied. Some examples where the splitting type as left module differed from the splitting type as right module
was constructed using explicit calculations. In this paper we give examples of sheaves of principal parts on projective space in any dimension where
the left $\O$-module structure differs from the right $\O$-module structure. Hence this gives the first examples that this phenomenon occur 
in any dimension. We also define Atiyah sequences and Atiyah classes of sheaves on arbitrary ringed topological spaces.
\end{abstract}

\maketitle

\tableofcontents

\section{Introduction} 

The aim of this paper is to give examples of sheaves of abelian groups equipped with two non-isomorphic structures as locally free sheaf
on projective space in any dimension. This was previously only known in dimension one (see \cite{maa100} and \cite{maa2}).
We do this using the first order sheaf of principal parts $\J^1(\O(l))$ of a line bundle $\O(l)$
where $l \geq $ is an integer (see Theorem \ref{mainsplit}). The proof relies on a well known result on the splitting type of the sheaf of prinicpal
parts as left module and the fact that the generalized Atiyah sequence is right split. We also show that the corresponding equivalence classes 
in algebraic  $\K$-theory are equal (see Theorem \ref{ktheory}). 

In the first section of the paper we give a definition of a generalized Atiyah sequence for an arbitrary ringed topological space $(X,\O)$
with respect to a derivation $d:\O\rightarrow \I$ where $\I$ is an abelianized sheaf of left and right $\O$-modules.
In the second section we prove that the left and right $\O$-module structure of the jet bundle of a line bundle on projective space in any dimension 
is non-isomorphic.
In the final section of the paper we consider the corresponding equvalence classes in algebraic $\K$-theory and prove these classes are equal.

\section{Generalized Atiyah sequences for ringed topological spaces}

The aim of this section is to define and prove basic properties of generalized Atiyah sequences for sheaves of modules on arbitrary 
ringed topological spaces.

Let in this section $X$ be an arbitrary topological space and let $\O$ be a sheaf of unital commutative rings on $X$.
Let $\I$ be an abelianized left and right $\O$-module. This means for every open set $U$ in $X$ and every section $a\in \O(U)$ and
$x\in \I(U)$ it follows $ax=xa$. Let $\E$ be a sheaf of left $\O$-modules. It follows $\E$ is canonically a sheaf of right $\O$-modules.
Let $d:\O\rightarrow \I$ be a derivation. This means $d$ is a morphism of sheaves of abelian groups satisfying the following formula on local sections 
$a,b\in \O(U)$: 
\[ d(ab)=ad(b)+d(a)b.\]

\begin{definition} \label{jet} Let $\J(\E)=\I\otimes_{\O}\E\oplus \E$ be the \emph{first order jet bundle} of $\E$.
\end{definition}

Let  $a\in \O(U)$ and $(x\otimes e,f)\in \J(\E)(U)$. Define the following structures of left and right $\O$-module on $\J(\E)$:
\begin{align}
&\label{j1} a(x\otimes e, f)=((ax)\otimes e+d(a)\otimes f, af)\\
&\label{j2}(x\otimes e, f)a=(x\otimes (ea), fa).
\end{align}

Define the following left and right $\O$-module structures on $\I\otimes_{\O}\E$:
\begin{align}
&\label{k1}a(x\otimes e)=(ax)\otimes e \\
&\label{k2}(x\otimes e)a=x\otimes(ea).
\end{align}

\begin{lemma} The definitions \ref{j1}, \ref{j2} define $\J(\E)$ as a left and right $\O$-module. Definition \ref{k1} and \ref{k2}
define $\I\otimes_{\O}\E$ as abelianized left and right $\O$-module.
\end{lemma}
\begin{proof} It is clear $\J(\E)$ and $\I\otimes_{\O}\E$ are left and right $\O$-modules. Assume $U$ is an open subset of $X$ and 
$a\in \O(U), x\otimes e\in \I\otimes_{\O}\E(U)$. We get
\[ a(x\otimes e)=(ax)\otimes e=(xa)\otimes e=x\otimes (ae)=x\otimes (ea)=(x\otimes e)a.\]
It follows $\I\otimes_{\O}\E$ is abelianized and the lemma is proved.
\end{proof}

We get an exact sequence
\begin{align}
&\label{atiyah}0\rightarrow \I\otimes_{\O}\E\rightarrow \J(\E)\rightarrow \E  \rightarrow 0
\end{align}
of sheaves of abelian groups on $X$.

\begin{theorem}\label{ringedtop}  The sequence \ref{atiyah} is an exact sequence of left and right $\O$-modules.  It is left split as sheaves of 
$\O$-modules if and only if there is an isomorphism $\J(\E)^{left}\cong \J(\E)^{right}$ of $\O$-modules.
\end{theorem}
\begin{proof} The proof is similar to the proof of Theorem 3.6 in \cite{maa10}.
\end{proof}

\begin{definition} Let $a(\E)\in \Ext^1_{\O}(\E,\I\otimes_{\O}\E)$ be the characteristic class defined by the sequence \ref{atiyah}.
It is called the \emph{generalized Atiyah class} of $\E$ with respect to the pair $(d,\I)$. Let the sequence \ref{atiyah} be the 
\emph{generalized Atiyah sequence} with respect to the pair $(d,\I)$.
\end{definition}

\begin{corollary} \label{corr} The following holds: $a(\E)=0$ if and only if $\J(\E)^{left}\cong \J(\E)^{right}$ as $\O$-modules.
\end{corollary}
\begin{proof} The proof is left to the reader as an exercise.
\end{proof}

Note: The definition of the Atiyah sequence was first done for holomorphic vector bundles on complex projective manifolds in Atiyah's paper 
\cite{atiyah}. In the paper \cite{karoubi} Karoubi defines the Atiyah sequence for commutative rings $A$ with respect to the pair $(d,\Omega^1_{A})$
where $d$ is the universal derivation. The definition given above where the Atiyah sequence is defined for any derivation $d:\O\rightarrow \I$
on any ringed topological space
gives a simultaneous generalization of a large class of Atiyah sequences: If $\I=\Omega^1_X(D)$ is the sheaf of logarithmic differentials where
$D\subseteq X$ is a normal crossing divisor and $d:\O_X \rightarrow \Omega^1_X(D)$ is the corresponding derivation it follows we get a 
\emph{logarithmic jet bundle} $\J(\E)(D)$ and a \emph{logarithmic Atiyah class} $a(\E)\in \Ext^1_{\O_X}(\E,\Omega^1_X(D)\otimes \E)$.
The logarithmic Atiyah class is the obstruction to the existence of a logarithmic connection
\[ \nabla: \E \rightarrow \Omega^1_X(D)\otimes \E.\]

\section{Examples: Jet bundles on projective space}

The aim of this section is to give new examples of sheaves of abelian groups on projective space in any dimension equipped with 
two non-isomorphic structures as locally free sheaf.

Let $K$ be an algebraically closed field of characteristic zero and let $V=K^{N+1}$.
Let $\P(V^*)$ be projective space on $V$.
Let $\O(1)$ be the tautological quotient bundle on $\P(V^*)$ and let $\O(-1)$ be its dual. 
Let $l\geq 1$ be an integer and let $\J(\O(l))$ be the first order jet bundle of $\O(l)$ as defined in \ref{jet}. It follows from \cite{maa10}
there is an isomorphism of left and right $\O_{\P(V^*)}$-modules
\[ \J(\O(l))\cong \J^1_{\P(V^*)}(\O(l)) \]
where $\J^1_{\P(V^*)}(\O(l))$ is the first order sheaf of principal parts of $\O(l)$.

\begin{lemma} \label{abel} The sheaf $\Omega^1_{\P(V^*)}\otimes \O(l)$ is an abelianized $\O_{\P(V^*)}$-module. The generalized 
Atiyah sequence \ref{atiyah} is right split as  $\O_{\P(V^*)}$-modules.
\end{lemma}
\begin{proof} Define the following map $s_U$ for any open subset $U$ in $X$:
\[ s_U: \O(l)|_U\rightarrow \J(\O(l))|_U \]
by
\[ s_U(x)=(0,x).\]
It follows
\[ s_U(xa)=(0,xa)=(0,x)a=s_U(x)a \]
hence $s_U$ is right $\O(U)$-linear for any open set $U$. One checks it splits the sequence \ref{atiyah} hence the first claim of 
the lemma is proved.
Assume $U=\Spec(A)\subseteq \P(V^*)$ is an open affine subset. It follows there is an isomorphism
\[ \Omega^1_{\P(V^*)}|_U\cong I/I^2 \]
of left and right $A$-modules, where $I\subseteq A\otimes_K A$ is the ideal of the diagonal. 
The $A$-module $I/I^2$ is generated by elements $\omega=db\otimes 1\otimes b-b\otimes 1$ where $b\in A$.
By definition the following holds:
\[ a\omega=a\otimes 1\omega \]
and
\[ \omega a=\omega 1\otimes a.\]
We get
\[ d(b)a-ad(b)=d(b)1\otimes a -a\otimes 1d(b)=d(b)(1\otimes a-a\otimes 1)=d(b)d(a)=0.\]
It follows 
\[d(b)a=ad(b) \]
in $I/I^2$.  Assume $L$ is a rank one locally free left $A$-module. The $A$-module $L$ has a canonical right $A$-module structure
defined by
\[ xa=ax.\]
We get for any $\omega \otimes x\in \Omega^1_A\otimes_A L$ and $a\in A$ the following calculation:
\[ a(\omega \otimes x)=(a\omega)\otimes x =(\omega a)\otimes x=\omega\otimes (ax)=\omega \otimes (xa)=(\omega \otimes x)a.\]
It follows $\Omega^1_A\otimes_A L$ is abelianized. The lemma follows.
\end{proof}

\begin{theorem} \label{mainsplit} Let $l\geq 1$ be an integer. There is no isomorphism
between $\J(\O(l))^{left}$ and $\J(\O(l))^{right}$ as $\O_{\P(V^*)}$-modules.
\end{theorem}
\begin{proof} By Lemma \ref{abel} the generalized Atiyah sequence is right split, hence it follows there is an isomorphism
\[ \J(\O(l))^{right}\cong \Omega^1_{\P(V^*)}\otimes \O(l)\oplus \O(l) \]
of right $\O_{\P(V^*)}$-modules.
From \cite{maa1} there is an isomorphism
\[ \J(\O(l))^{left}\cong \O(l-1)\oplus \cdots \oplus(l-1) \]
as left $\O_{\P(V^*)}$-modules. The theorem follows since $\O(l)$ cannot be a direct summand of $\O(l-1)\oplus \cdots \oplus \O(l-1)$.
\end{proof}

\begin{corollary} Let $l\geq 1$ be any integer. There is no connection
\[ \nabla:\O(l)\rightarrow \Omega^1_{\P(V^*)}\otimes \O(l) .\]
\end{corollary}
\begin{proof} The corollary follows from Theorem \ref{mainsplit} and Theorem 3.6 in \cite{maa10}.
\end{proof}

The sheaf of abelian groups $\J(\O(l))$ is locally free of rank $N+1$ as left and right $\O_{\P(V^*)}$-module on $\P(V^*)$ for any $N\geq 1$.
It follows the sheaf of abelian groups $\J(\O(l))$ is equipped with two non-isomorphic structures as 
locally free sheaf on projective space in any dimension. This was previously known to be true in dimension one (see \cite{maa100} and \cite{maa2}).

\section{Some classes in algebraic K-theory of projective space}

In the previous sections we gave examples of sheaves of abelian groups equipped with two non-isomorphic structures as locally free sheaf.
In this section we prove that the corresponding equivalence classes in algebraic $\K$-theory are equal.

Let $\P=\P(V^*)$ be projective space of dimension $N$ over $K$. Let $\K(\P)$ be the Grothendieck group of locally free finite rank sheaves
on $\P$.

\begin{proposition}\label{projectivespace} There is an isormophism
\[ \phi:\K(\P)\cong \Z\{1,t,t^2,\ldots , t^{N-1}\} \]
where
\[ \phi([\O(-d)])=\sum_{i=0}^{N-1} (-1)^i\binom{d}{i}t^i \]
and
\[ \phi([\O(d)])=\sum_{i=0}^{N-1} \binom{d+i-1}{i}t^i \]
where $t=1-[\O(-1)]$.
\end{proposition}
\begin{proof} By the projective bundle formula (see \cite{srinivas}, Proposition 5.18) there is an isomorphism
\[ \K(\P)\cong \K(\Spec(K))[h]/P(h) \]
where
\[ P(h)=h^N-[\wedge^1V]h^{N-1}+[\wedge^2 V]h^{N-2}+\cdots +(-1)^N[\wedge^N V]=(h-1)^N.\]
Here $h=[\O(-1)]$. We get
\[ \K(\P)\cong \Z[h]/(h-1)^N\cong \Z[h]/(1-h)^N\cong \]
\[ \Z[t]/t^N\cong \Z\{1,t,t^2,\ldots ,t^{N-1}\} .\]
We get
\[ [\O(-d)]=[\O(-1)]^d=(1+h-1)^d=\sum_{i=0}^{N-1}\binom{d}{i}(h-1)^i=\]
\[ \sum_{i=0}^{N-1}(-1)^i\binom{d}{i}(1-h)^i=\sum_{i=0}^{N-1}(-1)^i\binom{d}{i}t^i,\]
and the first claim is proved.
The following formula holds in general:
\[ \frac{1}{(1+a)^{k+1}}=\sum_{j=0}^{\infty} (-1)^j\binom{j+k}{j}a^j.\]
By definition $h^{-1}=[\O(1)]$ in $\K(\P)$. We get
\[ \phi([\O(d)])=\phi((h^{-1})^d)=\]
\[(\frac{1}{h})^d=\frac{1}{(1+h-1)^d}=\]
\[\sum_{j=0}^{N-1}(-1)^j\binom{j+d-1}{j}(h-1)^j=\sum_{j=0}^{N-1}\binom{d+j-1}{j}(1-h)^j=\]
\[ \sum_{j=0}^{N-1}\binom{d+j-1}{j}t^j.\]
and the proposition follows.
\end{proof}

Let $N=1$ and let $\E=\O(d_1)\oplus \cdots \oplus \O(d_e)$ be a locally free sheaf of rank $e$ on $\P^1$.
Let 
\[ deg(\E)=d_1+\cdots +d_e \]
the the \emph{degree} of $\E$.

\begin{corollary} \label{pline} There is an isomorphism
\[ \phi: \K(\P^1)\cong \Z\oplus \Z \]
defined by
\[ \phi([\E])=(deg(\E), rk(\E))\]
where $rk(\E)$ is the rank of the locally free sheaf $\E$.
\end{corollary}
\begin{proof} Let $N=1$ in proposition \ref{projectivespace}. We get
\[ \phi([\O(-d)])=\binom{d}{0}+\binom{d}{1}(h-1)=1-d(h-1)=1-dt.\]
We get similarly
\[ \phi([\O(d)])=1-d(h-1)=1+d(1-h)=1+dt.\]
It follows
\[\phi([\E])=\phi([\O(d_1)])+\cdots +\phi([\O(d_e)])=1+d_1t+\cdots +1+d_et=\]
\[ rk(\E)1+deg(\E)t.\]
The corollary follows.
\end{proof}

From \cite{maa2} the following holds:

\[ deg (\J^k(\O(l))^{left})=deg(\J^k(\O(l))^{right}) \]
and
\[ rk(\J^k(\O(l))^{left})=rk(\J^k(\O(l))^{right}) \]
for any integer $k\geq 1$. 

\begin{corollary} There is an equality
\[ [\J^k(\O(l))^{left}]=[\J^k(\O(l))^{right}] \]
of equivalence classes in $\K(\P^1)$.
\end{corollary}
\begin{proof} The corollary follows from the above discussion and Corollary \ref{pline}.
\end{proof}

Let in the following $X/S$ be a differentially smooth morphism of schemes. It follows for any finite rank locally free
$\O_X$-module we get fundamental exact sequences

\begin{align}
&\label{fundamental} 0\rightarrow \Sym^k(\Omega^1_{X/S})\otimes \E \rightarrow \J^k(\E) \rightarrow \J^{k-1}(\E) \rightarrow 0 
\end{align}

of left and right $\O_X$-modules generalizing the Atiyah sequence.
The sheaf of abelian groups $\Sym^k(\Omega^1_{X/S})\otimes \E$ has a canonical left and right $\O_X$-module structure.

\begin{lemma} \label{differentially} There is for any integer $k\geq 1$ an isomorphism
\[ (\Sym^k(\Omega^1_{X/S})\otimes \E)^{left}  \cong (\Sym^k(\Omega^1_{X/S})\otimes \E)^{right} \]
of $\O_X$-modules.
\end{lemma}
\begin{proof} A proof similar to the one in  Lemma \ref{abel} shows $\Omega^1_{X/S}$ is abelianized. It follows 
$\Sym^k(\Omega^1_{X/S})\otimes \E$ is abelianized for any integer $k\geq 1$.
\end{proof}

Let $\K(X)$ be the Grothendieck group of locally free finite rank sheaves on $X$.

\begin{theorem} \label{ktheory} There is an equality
\[ [\J^k(\E)^{left}]=[\J^k(\E)^{right}] \]
in $\K(X)$ for any integer $k\geq 1$.
\end{theorem}
\begin{proof} From Lemma \ref{differentially} and an induction we get the following:
\[ [\J^k(\E)^{left}]=\sum_{i=0}^k [(\Sym^k(\Omega^1_{X/S})\otimes \E)^{left} ]=\]
\[ \sum_{i=0}^k [(\Sym^k(\Omega^1_{X/S})\otimes \E)^{right} ] =[\J^k(\E)^{right}].\]
\end{proof}

We get examples of sheaves of abelian groups equipped with two non-isomorphic structures as locally free sheaf 
on projective space in any dimension. This is by Theorem \ref{ringedtop} 
detected by the generalized Atiyah class $a(\E)\in \Ext^1_{\O_X}(\E,\Omega^1_{X}\otimes \E)$.
The associated classes in algebraic $\K$-theory are equal.

Note: The exact sequence \ref{fundamental} only exist when $X/S$ is differentially smooth. In general there is an exact sequence
\[ 0\rightarrow \I^k/\I^{k+1}\otimes \E \rightarrow \J^k(\E) \rightarrow \J^{k-1}(\E)\rightarrow 0\]
where $\I\subseteq \O_{X\times X}$ is the ideal sheaf of the diagonal. When $X/S$ is differentially smooth it follows there is an
isomorphism
\[ \Sym^k(\Omega^1_X)\cong \I^k/\I^{k+1} \]
for every $k\geq 1$. The module $\I^k/\I^{k+1}$ is an abelianized sheaf of left and right $\O_X$-modules if and only if $k=1$.
Hence in the case when $X/S$ is not differentially smooth we do not neccessarily get the equality given in Theorem \ref{ktheory}.

\begin{example} Generalized Atiyah classes and Chern classes.\end{example}

Let $X$ be any scheme and $d:\O_X\rightarrow \Omega^1_X$ the universal derivation. Consider the generalized Atiyah sequence with 
respect to the pair $(d,\Omega^1_X)$. 

There is for any linebundle $\L$ on $X$ an isomorphism
\[ \Ext^1_{X}(\L,\Omega^1_{X}\otimes \L)\cong \Ext^1_{X}(\O,\Omega^1_{X}\otimes \L\otimes \L^*)\cong \]
\[ \Ext^1_{X}(\O,\Omega^1_{X} )\cong \H^1(X, \Omega^1_{X}).\]
There is a canonical map
\[ c_1: \Pic(X)\rightarrow \H^1(X, \Omega^1_{X}) \]
induced by the map $dlog:  \O_X^*\rightarrow \Omega^1_{X}$ 
and under this isomorphism we get an equality of characteristic classes $a(\L)=c_1(\L)$ in $\H^1(X, \Omega^1_{X})$

\begin{proposition} The following holds: $c_1(\L)=0$ if and only if there is an isomorphism
\[ \J(\L)^{left}\cong \J(\L)^{right} \]
of $\O_X$-modules.
\end{proposition}
\begin{proof} Since $a(\L)=c_1(\L)$ the proof follows from the discussion above and
Corollary \ref{corr}.
\end{proof}


\begin{thebibliography}{4}


\bibitem{atiyah} M. Atiyah, Complex analytic connections in fibre
  bundles, \emph{Trans. Amer. Math. Soc.} no. 85 (1957)

\bibitem{grothendieck} A. Grothendieck, EGA IV. \'{E}tude locale de sch\'{e}mas
  et des morphismes de sch\'{e}mas, I,\emph{Publ. Math. IHES} no. 20 (1964)
LMR0173675 (30:3885)

\bibitem{karoubi} M. Karoubi, Homologie cyclique et $K$-th\'{e}orie,
  \emph{Ast\'{e}risque} no 149 (1987)

\bibitem{maa1} H. Maakestad, A note on the principal parts on
  projective space and linear representations,
\emph{Proc. of the Amer. Math. Soc.} Vol. 133 no. 2 (2004) 

\bibitem{maa100} H. Maakestad, Modules of principal parts on the projective line, \emph{Arkiv for matematik} vol. 42, 2 (2004)

\bibitem{maa10} H. Maakestad, On jets, extensions and characteristic classes II, \emph{arXiv:1006.0593} (2010)


\bibitem{maa2} H. Maakestad, Principal parts on the projective line
  over arbitrary rings, \emph{Manuscripta Math.}, vol. 126, no. 4 (2008) 

\bibitem{srinivas} V. Srinivas, Algebraic K-theory, \emph{Progress in Mathematics, Birkhauser} no. 90 (1996) 


\end{thebibliography}
\end{document}